\documentclass[12pt]{article}
\usepackage{amsmath,amssymb,amsfonts,amsthm}
\usepackage{eucal} 

\newcommand{\lang}{\left\langle}
\newcommand{\rang}{\right\rangle}
\newcommand{\llang}{\lang\!\lang}
\newcommand{\rrang}{\rang\!\rang}

\newcommand{\lv}{\left |}
\newcommand{\rv}{\right |}
\newcommand{\zz}{{\mathfrak{z}}}
\newcommand{\vac}{v_\emptyset}
\newcommand{\bH}{\mathsf{H}}
\newcommand{\bS}{\mathsf{S}}
\newcommand{\bY}{\mathsf{Y}}

\newcommand{\bJ}{\mathsf{J}}

\newcommand{\cF}{\mathcal{F}}
\newcommand{\bw}{\mathsf{w}}
\newcommand{\bZ}{\mathsf{Z}}

\newcommand{\bGg}{\mathsf{G}_{\textup{\tiny\sc GW}}}
\newcommand{\bGd}{\mathsf{G}_{\textup{\tiny\sc DT}}}

\newcommand{\bn}{\mathsf{n}}

\newcommand{\bGm}{\mathbf{\Gamma}}

\newcommand{\MM}{\mathsf{M}}

\newcommand{\C}{\mathbb{C}}
\newcommand{\Q}{\mathbb{Q}}
\newcommand{\Z}{\mathbb{Z}}
\newcommand{\cO}{\mathcal{O}}

\newcommand{\Pp}{\mathbf{P}^1}

\DeclareMathOperator{\Res}{Res}

\DeclareMathOperator{\sgn}{sgn}

\DeclareMathOperator{\Aut}{Aut}

\DeclareMathOperator{\End}{End}

\newtheorem{Theorem}{Theorem}
\newtheorem{Lemma}{Lemma}
\newtheorem{Corollary}[Lemma]{Corollary}
\newtheorem{Proposition}[Lemma]{Proposition}

\newcommand{\bP}{\mathsf{P}}

\newcommand{\bff}{\mathbf{f}}
\newcommand{\yb}{\mathbf{y}}
\newcommand{\jb}{\mathbf{j}}
\newcommand{\mb}{\mathbf{m}}
\newcommand{\lc}{<_c}
\newcommand{\gc}{>_c}

\begin{document}
\title{The quantum differential equation of
 the Hilbert scheme of 
points in the plane}
\author{A.~Okounkov and R.~Pandharipande} 
\date{} \maketitle


\section{Introduction}

\subsection{Overview}

In the study of the quantum cohomology of the Hilbert scheme of 
points in the plane \cite{QCHS}, as well as in the 
Gromov-Witten/Donaldson-Thomas theories of threefolds \cite{jbrp,moop,GWDT}, 
 certain linear ODEs with remarkable 
properties arise naturally. 
These ODEs generalize the Schr\"odinger 
equation for the quantum Calogero-Sutherland operator and 
are the focus of the present paper. Following the tradition 
in their field of origin, we call them quantum differential 
equations, or QDEs for short, even though there is nothing 
quantum per se in these differential equations. They are
linear ODEs with regular singularities, very classical 
objects indeed. 

Two special values of the independent variable $q$ play 
a special role. These are $q=-1$ and $q=0$, which may be 
called the Gromov-Witten and Donaldson-Thomas points, 
respectively. The point $q=-1$ is nonsingular. We prove 
that the monodromy based at $q=-1$ is a polynomial in 
$e^{2\pi i t_i}$, where $t_1$ and $t_2$ are the parameters
of the QDE. 

Next we solve the connection problem for 
the points $q=-1$ and $q=0$. The point $q=0$ is singular 
with the residue being essentially the Calogero-Sutherland
operator. Its eigenfunctions are known as the Jack 
symmetric functions. We prove that transported to $q=-1$ 
by the QDE, these become, up to normalization, 
Macdonald polynomials with parameters $e^{2\pi i t_i}$. 

While both the objects and the results of this paper
belong to the world of combinatorics and differential 
equations, our proofs require geometric input at several 
key points. Perhaps a deeper understanding of the integrable
structures underlying the QDEs will lead to more direct proofs. 

The appearance of Macdonald polynomials in the connection 
problem strongly suggests a relation to the equivariant
$K$-theory of the Hilbert scheme. This relation is further
pursued in \cite{BO}.

\subsection{Fock space and symmetric functions}
\label{fs}

The most natural way to write down our ODEs is in terms of 
creation and annihilation operators acting on the Fock 
space. 

We follow the conventions of \cite{QCHS}. By definition, the Fock space $\cF$ 
is freely generated over $\C$ by commuting 
creation operators $\alpha_{-k}$, $k\in\Z_{>0}$,
acting on the vacuum vector $\vac$. The annihilation 
operators $\alpha_{k}$, $k\in\Z_{>0}$, kill the vacuum 
$$
\alpha_k \cdot \vac =0,\quad k>0 \,,
$$
and satisfy the commutation relations
$$
\left[\alpha_k,\alpha_l\right] = k \, \delta_{k+l}\,. 
$$
A natural basis of $\cF$ is given by 
the vectors  
\begin{equation}
  \label{basis}
  \lv \mu \rang = \frac{1}{\zz(\mu)} \, \prod \alpha_{-\mu_i} \, \vac \,.
\end{equation}
indexed by partitions 
$\mu$. Here, $$\zz(\mu)=|\Aut(\mu)| \, \prod \mu_i$$ is the usual 
normalization factor. 

The linear map 
\begin{equation}
p_\mu \mapsto  \zz(\mu) \, \lv \mu \rang \,,\label{symmcF}
\end{equation}
where 
$$
p_\mu = \prod_k \sum_{i=1}^\infty z_i^{\mu_k} \,,
$$
identifies $\cF$ with symmetric functions of the variables
$$z_1,z_2,z_3,\dots\, .$$ 
Symmetric 
functions of fixed degree form eigenspaces of 
the {\em  energy operator}: 
$$
|\cdot| = \sum_{k>0} \alpha_{-k} \, \alpha_k \,.
$$

\subsection{The QDE}

The central object of this paper is the differential equation 
\begin{equation}
  \label{de1}
   q\frac{d}{dq} \, \Psi = \MM \, \Psi\,, \quad \Psi\in \cF\,,
\end{equation}
where the operator $\MM$ is given by 
\begin{multline}
  \label{theM} 
\MM(q,t_1,t_2) = (t_1+t_2) \sum_{k>0} \frac{k}{2} \frac{(-q)^k+1}{(-q)^k-1} \,
 \alpha_{-k} \, \alpha_k  + \\
\frac12 \sum_{k,l>0} 
\Big[t_1 t_2 \, \alpha_{k+l} \, \alpha_{-k} \, \alpha_{-l} -
 \alpha_{-k-l}\,  \alpha_{k} \, \alpha_{l} \Big] \,.
\end{multline}
The variables $t_1$ and $t_2$ in
\eqref{theM} are parameters. 
Note that the $q$-dependence of $\MM$ is only in the first sum in \eqref{theM}
which acts diagonally in the basis \eqref{basis}. 
The two terms in the 
second sum in \eqref{theM} are known respectively as the splitting and joining terms. 

The operator 
  $\MM$ commutes with the energy operator 
$|\cdot|$, therefore the equation splits into a direct sum of 
finite-dimensional ODEs. It is more convenient to study 
the equivalent equation 
\begin{equation}
  \label{de2}
   q\frac{d}{dq} \, \Psi = \MM_D \Psi\,, \quad \Psi\in \cF\,,
\end{equation}
where
  \begin{equation}
    \label{MM_D}
    \MM_D = \MM - \frac{t_1+t_2}{2} \, \frac{(-q)+1}{(-q)-1} \, 
|\,\cdot\,| \,\,.
  \end{equation}
The equation \eqref{de2} is the quantum differential equation 
for the Hilbert scheme of points, see \cite{QCHS}. One advantage of
\eqref{de2} is that $q=-1$ is a regular point for \eqref{de2}. 

In this paper, we normalize everything exactly as in \cite{QCHS} 
even though some simplifications such as replacing $q$ by 
$-q$ in \eqref{theM} may seem obvious.

\subsection{Calogero-Sutherland operator}

The quantum-mechanical Calogero-Sutherland
operator, 
\begin{equation}
  \label{HCS}
  \bH_{CS}= \frac12 \sum_i 
\left( z_i \frac{\partial}{\partial z_i}\right)^2 +
\theta(\theta-1) \sum_{i<j} \frac1{|z_i-z_j|^2} \,,
\end{equation}
describes  particles moving on the torus 
$|z_i|=1$ interacting via the potentials $|z_i-z_j|^{-2}$. The parameter
$\theta$ adjusts the strength of the interaction. 
The function 
$$
\phi(z) = \prod_{i<j} (z_i-z_j)^\theta
$$
is an eigenfunction of $\bH_{CS}$, and the operator
$\phi \, \bH_{CS} \, \phi^{-1}$ preserves the space of symmetric
polynomials in the variables $z_i$. Therefore,
via the identification \eqref{symmcF}, 
the operator
$\phi \, \bH_{CS} \, \phi^{-1}$
acts on the Fock space.

A direct computation shows  the operator  $\phi \, \bH_{CS} \, \phi^{-1}$ equals
\begin{equation}
 \Delta_{CS}= \frac{1-\theta}{2}\sum_k k \, \alpha_{-k} \, \alpha_k  + \\ 
\frac12 \sum_{k,l>0} 
\Big[\alpha_{-k-l} \, \alpha_{k} \, \alpha_{l} + \theta \, 
 \alpha_{k+l}\,  \alpha_{-k} \, \alpha_{-l} \Big] \, 
\end{equation}
modulo scalars and a multiple of the momentum 
operator $\sum_i z_i \frac{\partial}{\partial z_i}$, see \cite{St}.
We find
\begin{equation}
\label{MDCS}
  \MM(0) = 
- \,t_1^{\ell(\,\cdot\,)+1} \, 
\Delta_{CS}\Big|_{\theta=-t_2/t_1} \, t_1^{-\ell(\,\cdot\,)}  \,, 
\end{equation}
where $\ell(\mu)$ is the  number of 
parts of the partition $\mu$ and  
$\ell(\,\cdot\,)$ is the diagonal operator with eigenvalues
$\ell(\mu)$ in the basis $\lv \mu \rang$.

The well-known duality $\theta\mapsto1/\theta$
in the Calogero-Sutherland model corresponds to 
the permutation of $t_1$ and $t_2$. 

Formula \eqref{MDCS} implies the behavior of \eqref{de1} near the 
regular singularities $q=0,\infty$ is described by the Schr\"odinger
equation for \eqref{HCS}. The connection problem for these two 
points may be viewed as a scattering produced by nonstationary 
terms in \eqref{de1}. This problem will be considered in Section \ref{s_scatt}.

\subsection{An application}

The solution of the connection problem may be combined with the results 
of \cite{moop} to give a box-counting formula for triple Hodge integrals. 
While we will not
reproduce the required  formulas here, the essence summarized 
as follows. 

In \cite{moop} a certain enumerative object, called the \emph{capped 
vertex} was introduced and shown to satisfy the GW/DT correspondence.
Its ingredients on the GW and DT sides are quite different. On the 
GW side, it encapsulates the general triple Hodge integrals, 
repackaged using fundamental solutions of our QDE, normalized at $q=-1$. 

On the DT side, the expression is essentially combinatorial, 
its main ingredient being a certain weighted count of 3-dimensional 
partitions known as the \emph{equivariant vertex}. It is similarly
decorated by the fundamental solutions of the QDE, but this time 
normalized at $q=0$. Therefore, it is precisely the connection 
formulas for the QDE that relate the two expressions. 

\subsection{Acknowledgments} 
We thank R.~Bezrukavnikov, P.~Etingof, N.~Katz, D.~Maulik,  N.~Nekrasov, and A.~Oblomkov for  
valuable discussions. The results presented here date back several years. See, in particular, 
\cite{BO} for further developments. 

Both authors were partially 
supported by the Packard 
foundation and the NSF.

\section{Basic properties}

\subsection{Singularities}
For $\Psi\in\cF$ of energy $n$, the equation \eqref{de2} is 
 a linear first order ODE in $p(n)$ unknowns, where $p(n)$
is the number of partitions of $n$. It has regular singularities. 
These are $q=0,\infty$, and solutions $\zeta$ of 
$$
(-\zeta)^m = 1, \quad m=2,\dots,n  \,,
$$
excluding $q=-1$. 

For example, for $n=3$, in the basis \eqref{basis} 
ordered lexicographically, the matrix $\MM_D$ takes the 
form
$$
\begin{bmatrix}
3(t_1+t_2) 
\dfrac{q^2-1}{q^2-q+1} &-3&0\\
2 t_1 t_2 & (t_1+t_2) \dfrac{q+1}{q-1} &-1 \\
0&3 t_1 t_2 &0
\end{bmatrix} \,.
$$

\subsection{Residues at $q=0$}

Equation \eqref{MDCS} relates the residue of \eqref{de2} at $q=0$ 
to the Calogero-Sutherland operator. In particular, the 
eigenfunctions of $\MM(0)$ are, up to normalization, Jack 
symmetric functions. 

More precisely, let $J_\lambda\in \cF$ be integral form of the Jack 
symmetric function depending on the parameter $\alpha=1/\theta$ as in \cite{Mac}. We define
\begin{equation}
\bJ^\lambda = 
t_2^{|\lambda|} \, t_1^{\ell(\cdot)} \, J_\lambda\big|_{\alpha = - t_1/t_2} \,.\label{dbJ}
\end{equation}
These are eigenfunctions of $\MM_D(0)$ with eigenvalues 
\begin{equation}
  \label{def_c}
 - c(\lambda;t_1,t_2) = - \sum_{(i,j)\in \lambda} \Big[
(j-1) t_1 + (i-1) t_2 \Big] \,. 
\end{equation}
The normalization is such that 
\begin{equation}
  \label{normJ}
  \bJ^\lambda = |\lambda| ! \, (t_1 t_2)^{|\lambda|}
\, \lv 1^{|\lambda|} \rang + \dots \,.
\end{equation}
For example, we have
\begin{equation*}
\bJ^{(k)} = k! \, t_1^k\,  \sum_{|\mu|=k} (-1)^{k-\ell(\mu)}\, t_2^{\ell(\mu)} \, 
\lv \mu \rang \,.
\end{equation*}
In general, the coefficient of $\lv \mu \rang$ in the 
expansion of $\bJ^\lambda$ is $(t_1 t_2)^{\ell(\mu)}$ 
times a polynomial in $t_1$ and
$t_2$ of degree $|\lambda|-\ell(\mu)$.

Geometrically, $\bJ^\lambda$ correspond to classes of monomial ideals in
the equivariant cohomology of the Hilbert scheme of points, see 
 \cite{vass,LQW}. From 
this point of view, the symmetry 
\begin{equation}
\bJ^\lambda(t_2,t_1)=\bJ^{\lambda'}(t_1,t_2)\label{symmJ},
\end{equation}
where $\lambda'$ is the transpose of $\lambda$,
is obvious.

By the general theory of ODEs, we can construct a solution
$$
\Psi = \bY^\lambda(q) \, q^{-c(\lambda;t_1,t_2)}\,, \quad \bY^\lambda(q) \in \C[[q]] 
$$
of \eqref{de2} which converges for $|q|<1$ and satisfies 
$$
\bY^\lambda(0) = \bJ^\lambda \,.
$$
The symmetry \eqref{symmJ} generalizes to 
$$
\bY^\lambda(t_2,t_1)=\bY^{\lambda'}(t_1,t_2)\,.
$$

\subsection{Other residues}

We find
\begin{equation}
  \label{q->1/q}
  \MM_D(q^{-1}) = - (-1)^{\ell(\,\cdot\,)}\, \MM_D(q) \,  (-1)^{\ell(\,\cdot\,)}\,.
\end{equation}
Therefore, the exponents at $q=\infty$
are the same as exponents at $q=0$.

If $\zeta\ne -1$ is a root of unity, then 
\begin{equation}
\Res_{q=\zeta} q^{-1} \, \MM_D(q) = (t_1+t_2) 
\sum_{\{k|(-\zeta)^k=1\}} \alpha_{-k} \, \alpha_k \,,\label{res_root}
\end{equation}
which is diagonal in the basis \eqref{basis}. 
In particular, the exponents at the 
roots of unity are positive integer multiples of 
the number $(t_1+t_2)$, which we call the  
\emph{level} of the equation. 

For positive integer levels, the solutions will be seen to be regular
at the roots of unity, implying that the $\bY^\lambda$ 
are polynomial.

\subsection{Unitarity}

The standard inner product on cohomology  induces a certain 
bilinear inner product on $\cF$ with respect to which $\MM$ is symmetric. 

As will be explained below, it is more convenient for us to work with 
a certain Hermitian inner product --- an inner product 
\emph{antilinear} in the second entry. This means that 
$$
\lang a f, g\rang = a\lang f, g\rang\,, \quad a\in\C(t_1,t_2)
$$
and
$$
\lang f, g\rang = \overline{\lang g, f\rang} \,,
$$
where, by definition
$$
\overline{a(t_1,t_2)} = a(-t_1,-t_2) \,.
$$
Specifically, we consider the Hermitian product defined on basis 
vectors by 
\begin{equation}
  \label{inner_prod}
  \lang \mu | \nu \rang = 
\frac{1}{(t_1 t_2)^{\ell(\mu)}} 
\frac{\delta_{\mu\nu}}{\zz(\mu)} \,. 
\end{equation}
We then have 
\begin{equation}
  \label{adjoint}
  \left(\alpha_{k}\right)^* = (t_1 t_2)^{\sgn(k)} \, 
\alpha_{-k}\,,
\end{equation}
and hence $\MM$ is skew-Hermitian: 
\begin{equation}
  \label{Madj}
  \MM^* = -\MM \,\, .
\end{equation}
The following is an immediate corollary

\begin{Proposition}\label{prop_U} The connection 
$$
\nabla(t_1,t_2) = q\frac{d}{dq} - \MM_D(q;t_1,t_2) 
$$
is unitary. 
\end{Proposition}

The polynomials $\bJ^\lambda$ are orthogonal with 
respect to \eqref{inner_prod}. Their norm is given by the product of the 
tangent weights:  
\begin{equation}
  \label{normJJ}
\|\bJ^\lambda\|^2 = \prod_{\textup{tangent weights $\bw$}} \bw
\end{equation}
We recall that the tangent weights $\bw$ to the Hilbert scheme at the monomial ideal 
indexed by $\lambda$ are given by 
\begin{equation}
\{\bw\} = \big\{t_1(a(\square)+1) - t_2 l(\square), -t_1 a(\square) + t_2 (l(\square)+1)\big\}_{\square\in \lambda}\,,
\label{tang_weights}
\end{equation}
where $a(\square)$ and $l(\square)$ denote the arm-length $\lambda_i-j$ and 
leg-length $\lambda'_j-i$ of a square $\square=(i,j)$ in a the diagram $\lambda$. 
Proposition \eqref{prop_U} implies the following 

\begin{Corollary}
$$
\lang \bY^\lambda(q),\bY^\mu(q)\rang =  \delta_{\lambda\mu} \, \|\bJ^\lambda\|^2 \,.
$$  
\end{Corollary}

Proposition \ref{prop_U} may also be phrased as a relation 
between the fundamental solutions on opposite levels. 

There is a geometric reason why the Hermitian inner product \eqref{inner_prod} 
is preferred. As explained in \cite{moop}, it is natural do define inner products on 
$\cF$ as equivariant GW/DT invariants of $\Pp\times\C^2$ relative to
the fibers over 
$0,\infty\in\Pp$. In that setting, our parameters $t_1$ and $t_2$ correspond to 
the ratios of torus weights in $\C^2$-direction to the weight in the $\Pp$-direction. 
But, clearly, the tangent spaces to $\Pp$ at $0$ and $\infty$ have \emph{opposite} 
torus weights.

\section{Monodromy}

\subsection{Intertwiners}\label{s_monod}

\subsubsection{}

A crucial role in what follows will be played by certain operators 
$$
\bS(a,b)\in\End(\cF)\otimes \Q(q,t_1,t_2)\,, \quad (a,b)\in \Z^2\,, 
$$ 
that intertwine the monodromy
of $\nabla(t_1,t_2)$ and $\nabla(t_1-a,t_2-b)$, that is,  
satisfy  
\begin{equation}
  \label{intertw}
\nabla(t_1,t_2) \, \bS(a,b)  = \bS(a,b) \, \nabla(t_1-a,t_2-b)\,.
\end{equation}

The construction of these operators is geometric and is provided 
by the Gromov-Witten and Donaldson-Thomas theories of local curves. 
The equivalence of these two theories, a very special case of general GW/DT 
conjectures of \cite{mnop}, was proven in \cite{GWDT}. 
Specifically, the GW/DT partition function of the total space of
the $\cO(a)\oplus\cO(b)$ bundle over $\Pp$, relative the fibers 
over $0,\infty\in\Pp$ defines an operator on $\cF$. Up to 
normalization, these are the intertwiners $\bS$. 

A formula for the intertwiners $\bS$ in terms of the fundamental
solution of QDE may be derived using equivariant localization in 
GW and DT theories, respectively. We refer to \cite{GWDT} for technical 
details and state only the final result here. 

\subsubsection{}

Let $\Phi(q;t_1,t_2)$ be the fundamental solution of the QDE
\begin{equation}
  \label{qfund}
  q\frac{d}{dq} \, \Phi(q;t_1,t_2) = \MM_D \, \Phi(q;t_1,t_2)\,, 
\end{equation}
normalized by
$$
\Phi(-1;t_1,t_2)=1 \,.
$$
Denote 
$$
g(x,t) = \frac{x^{tx}}{\Gamma(tx)}\,, \quad x>0 \,.
$$
and let $\bGg$ denote the diagonal operator with the following eigenvalues
\begin{equation}
  \label{glue}
  \bGg(t_1,t_2) \, \lv \mu \rang =
\prod_{i} g(\mu_i,t_1)\, g(\mu_i,t_2)
\, \lv \mu \rang \,.
\end{equation}

\subsubsection{}

Let $\bY$ denote the matrix formed by the vectors $\bY^\lambda$. Let $\bGd(t_1,t_2)$ 
be the diagonal matrix with eigenvalues 
$$
q^{-c(\lambda,t_1,t_2)} \prod_{\textup{tangent weights $\bw$}} \frac{1}{\Gamma(\bw+1)}\,,
$$
where $c(\lambda,t_1,t_2)$ is the sum of $(t_1,t_2)$-contents defined in \eqref{def_c} 
and the product ranges over the tangent weights $\bw$ to the Hilbert scheme at
the monomial ideal labeled by the partition $\lambda$ as in \eqref{tang_weights}.  
We fix a branch of the multivalued functions $q^{-c(\lambda,t_1,t_2)}$ by 
\begin{equation}
  q^{-c(\lambda,t_1,t_2)}\Big|_{q=-1} = e^{\pi i c(\lambda,t_1,t_2)} \label{branch}
\end{equation}
The matrix $\bY \, \bGd$ is a solution of the
QDE with a particular normalization near $q=0$. 

\subsubsection{}
Equivariant localization of relative GW/DT invariant of $\cO(a)\oplus\cO(b)$ yields the 
following result

\begin{Theorem} For $(a,b)\in\Z^2$, there exists $\bS(a,b)\in\End(\cF)\otimes \Q(q,t_1,t_2)$ 
such that 
\begin{align}
\bS(a,b)&= \Phi(t_1,t_2)\,  \bGg(t_1,t_2)\,  \bGg(t_1-a,t_2-b)^{-1} \, \Phi(t_1-a,t_2-b)^{-1} \notag \\
        &= \bY(t_1,t_2) \,  \bGd(t_1,t_2)\,  \bGd(t_1-a,t_2-b)^{-1} \, \bY(t_1-a,t_2-b)^{-1} \label{Sform}\,. 
\end{align}
\end{Theorem}

\subsubsection{}

Formulas \eqref{Sform} are derived as follows. Let $\bZ'_\textup{\tiny \sc DT}(\lambda,\mu)$ 
denote the reduced Donaldson-Thomas partition function of the total space 
of  
$$\cO(a)\oplus\cO(b) \rightarrow \Pp,$$
 relative to the fibers over $0,\infty\in\Pp$. The partitions 
$\lambda$ and $\mu$ record the tangency to these fibers. The degree $d=|\lambda|=|\mu|$
of the curve is implicit in the notation. By definition $\bZ'_\textup{\tiny \sc DT}(\lambda,\mu)$
is a generating function counting 1-dimensional 
subschemes $C$ with weight $q^{\chi(\cO_C)}$. 
By Theorem 3 of \cite{GWDT},
\begin{equation}
  (-q)^{-d(2+a+b)/2} \, \bZ'_\textup{\tiny \sc DT}(\lambda,\mu) =
(-iu)^{d(a+b)+\ell(\lambda)+\ell(\mu)} \, \bZ'_\textup{\tiny \sc GW}(\lambda,\mu)
\label{GW=DT}
\end{equation}
after the change of variables $q=-e^{iu}$. Here,
 $\bZ'_\textup{\tiny \sc GW}(\lambda,\mu)$ 
is the generating function for the Gromov-Witten counts
weighting 
genus $g$ curves by
$u^{2g-2}$. 
Theorem 2 of \cite{GWDT} shows $\bZ'_\textup{\tiny \sc DT}(\lambda,\mu)$ is 
a rational function of $q$. 

The equality \eqref{GW=DT} is proven in \cite{GWDT} in equivariant cohomology 
with respect to the fiberwise  $(\C^*)^2$-action. A stronger result that involves 
the full $(\C^*)^3$ of automorphisms is proven in Proposition 1 of \cite{moop}. 

Both sides of \eqref{GW=DT} may be computed by equivariant localization, see 
for example, Section 2 of \cite{moop} for a summary. The
formula involves the
\emph{rubber integrals} as well as certain \emph{edge weights}. The relation 
of rubber integrals to the fundamental solution of our QDE, and hence to 
the matrix $\bY$, is discussed in Section 11.2 of \cite{GWDT}. (Note a scalar 
difference between operators $\MM$ and $\MM_D$.) 
Similarly, Gromov-Witten rubber 
integrals lead to the matrix $\Psi$. 

The Donaldson-Thomas edge weights are described in Section 4 of \cite{mnop}. 
In both theories, edge weights are rational function of equivariant parameters that factor into 
linear factors. They may be rewritten as ratios of $\Gamma$-functions, as
they appear in the middle of \eqref{GW=DT}. After making all terms explicit,
cancelling a common prefactor, and scaling the equivariant weight of 
$T_0\Pp$ to equal $1$, we arrive at the claim of the Theorem.

\subsubsection{}

The intertwining property \eqref{intertw} is obvious from \eqref{Sform} as are the following 
composition property 
$$
\bS(a_1,b_1;t_1,t_2) \, \bS(a_2,b_2; t_1-a_1,t_2-a_2) = \bS (a_1+a_2,b_1+b_2; t_1,t_2) 
$$
and symmetry
$$
\bS(a,b;t_1,t_2) = \bS(b,a; t_2,t_1) \,.
$$
We further have 

\begin{Proposition}
The intertwiner $\bS(a,b)$ is a Laurent polynomial in $q$ when $a+b\ge 0$. 
\end{Proposition}

\begin{proof}
The matrix $\bS(a,b)$ is a rational solution of an ODE \eqref{intertw} with regular singularities. 
Its possible singularities are controlled by the integer exponents of \eqref{intertw}. By 
\eqref{res_root} all integer exponents at roots of unity are nonnegative when $a+b\ge 0$. 
\end{proof}

\subsection{}

Since the matrix $\Phi$ is nonsingular away from the 
singularities of the differential equation, 
the additional singularities of $\bS$ are determined by the 
diagonal matrix $\bGg$. 
We conclude the following result.

\begin{Theorem}\label{T_monodr} 
The connections $\nabla(t_1,t_2)$ and $\nabla(t_1,t_2-1)$ have 
isomorphic monodromy provided
$$
t_2\ne \frac{r}{s}\,, \quad 0< r\le s \le n \,.
$$
\end{Theorem}

At zero level, $t_1+t_2=0$, the matrix $\MM_D$ is constant in $q$ 
and, hence, the monodromy is abelian. 
Furthermore, it is semisimple unless $t_1=t_2=0$. 
We can use Theorem \ref{T_monodr} repeatedly to get
from an integer level to zero level.

\begin{Corollary}
At integer level $t_1+t_2=l\in \Z$, the monodromy is 
abelian. Additionally, the monodromy is semisimple provided 
\begin{equation}
  \label{Ab_sem}
  t_1 \ne r/s
\end{equation}
with $0<s\le n$ and 
$$
r =
\begin{cases}
 1,\dots,ls-1\,, \quad & l>0 \,, \\
ls,\dots,0\,, \quad & l\le 0 \,. 
\end{cases}
$$
\end{Corollary}

\begin{proof}
By Theorem \ref{T_monodr}, we need only prove the
monodromy is semisimple for 
$$
(t_1,t_2)=(l,0)\,, \quad l>0 \,.
$$ 
For the rest of this proof we use $>$ to denote the 
following partial ordering on partitions of a fixed number $n$. It is the 
transitive closure of the following relation: $\lambda>\mu$ 
if $\lambda_i=\mu_j+\mu_k$ for some $i,j,k$, all other parts 
being identical in the two partitions. 
The matrix $\MM_D(q;l,0)$ is upper-triangular in the 
basis \eqref{basis} with respect to this partial ordering. 

We look for solutions of the form 
\begin{equation}
  \label{form_sol}
    B(q)\, q^{\MM_D(0;l,0)} \,,
\end{equation}
where $B(q)$ is holomorphic in $q$ and upper-triangular with $1$'s on the diagonal. The
matrix $\MM_D(0;l,0)$ is semisimple with 
eigenvalues of the form 
$$
-c(\lambda;l,0)= - l \sum_{i} \frac{\lambda_i^2}{2} \,.
$$
If $\lambda > \mu$ and $l>0$ then, 
$$
c(\lambda;l,0) - c(\mu;l,0) > 0 
$$
because $(a+b)^2>a^2+b^2$ for $a,b>0$. 
Thus, the eigenvalues of $\MM_D(0;l,0)$ are 
strictly increasing down the diagonal. Hence, by a straightforward argument,
a formal solution \eqref{form_sol}
can always be calculated order by order in $q$. 
\end{proof}

\subsection{}

Consider the operator
$$
\bGm \, \lv \mu \rang  = \frac{(2\pi i)^{\ell(\mu)}}{\prod \mu_i} \, \bGg(t_1,t_2) \lv \mu \rang 
$$
and the connection
$$
\nabla^\Gamma = \bGm \, \nabla \, \bGm^{-1} \,.
$$
The monodromy of this connection with base point $q=-1$ defines a representation 
\begin{equation}
  \pi_1(\C\setminus \{\textup{singularities}\},-1) \to \Aut(\cF)
\label{monodromy_rep}
\end{equation}

\begin{Theorem}\label{T_polynom} 
Matrix elements of \eqref{monodromy_rep} are Laurent polynomials in 
$$
T_i = e^{2\pi i t_i}\,, \quad i=1,2 \,.
$$  
Further, the monodromy \eqref{monodromy_rep} 
is unitary with respect to the Hermitian form defined by
\begin{equation}
  \label{Herm_prod}
  \llang \mu | \nu \rrang = 
\frac{\delta_{\mu\nu}}{\zz(\mu)} 
\, 
\prod_i \left(T_1^{\mu_i/2} - T_1^{-\mu_i/2}\right)\left(T_2^{\mu_i/2} - T_2^{-\mu_i/2}\right) \,.
\end{equation}
\end{Theorem}

\noindent 
The Hermitian form above is anti-linear with respect to the 
involution 
$$
\overline{T_i}=T_i^{-1}\,, \quad i=1,2 \,.
$$
We note that the Hermitian product \eqref{Herm_prod} is essentially the Macdonald inner 
product as well as the natural inner product in the K-theory of the Hilbert scheme of points, see \cite{Hai}. 
This will be revisited below and, more fully, in \cite{BO}.

\begin{proof}
Let $\gamma(t_1,t_2)$ be the monodromy of $\nabla$ along a loop based at $q=-1$. 
By \eqref{intertw}, we have 
$$
\gamma(t_1-a,t_2-b) = S(-1;a,b) \, \gamma(t_1,t_2) S(-1;a,b)^{-1} \,.
$$
On the other hand, 
$$
S(-1;a,b) =\bGm(t_1-a,t_2-b)^{-1} \,  \bGm(t_1,t_2)  \,.
$$
Hence $\bGm\, \gamma \, \bGm^{-1}$ is invariant under $t_i\mapsto t_i+1$ and, 
therefore, is a meromorphic function of $T_1$ and $T_2$. 

Since $\nabla$ depends polynomially on $t_1$ and $t_2$, its monodromy is 
an entire function of $t_1$ and $t_2$. The matrix $\bGm$ is holomorphically 
invertible for $\Re t_i > 0$, hence the monodromy of $\nabla^\Gamma$ is holomorphic 
there. By periodicity, it is holomorphic for all $T_1,T_2\in\C^*$.  

The solutions of $\nabla$ grow at most exponentially as $\Im t_i\to\infty$. {F}rom the 
Stirling formula, we have 
\begin{equation}
\ln | \Gamma(x+iy)| = - \frac{\pi}{2} |y| + O(\ln |y|) 
\,,\quad 
y\to \pm\infty \,, \label{Stirling}
\end{equation}
hence the monodromy of $\nabla^\Gamma$ grows at most polynomially as $T_i\to 0,\infty$. 
Therefore, $\bGm\, \gamma \, \bGm^{-1}$ is a Laurent polynomial in the $T_i$'s. 

By construction, the monodromy is unitary with respect to the new Hermitian form defined by
$$
\llang \mu | \nu \rrang = \lang \mu \lv \bGm(-t_1,-t_2) \bGm(t_1,t_2) \rv \nu \rang 
$$
Using the formula
$$
- \frac{2 \pi i}{\Gamma(x)\Gamma(-x)} = x \, (e^{\pi i x} - e^{-\pi i x}) 
$$
we obtain \eqref{Herm_prod}. 
\end{proof}

\subsection{$\Gamma$-factors}

This section is full of $\Gamma$-factors and one may wonder what is their deeper 
meaning. Technically, they arise as edge-weights in GW and DT localization 
formulas. Further, associated to the intertwiner operator $\bS(a,b)$ is a
system of \emph{difference equations} 
\begin{equation}
  \bS(a,b) \Psi(t_1 - a,t_2 -b ) = \Psi(t_1,t_2)\label{Seq}
\end{equation}
which is compatible with the differential equation \eqref{de1}. The variable 
$q$ of the original differential equation \eqref{de1} is a parameter in 
the difference equation \eqref{Seq}. The points $q=0,-1$ are special 
for this difference equation in that for those values of the parameter
it becomes abelian and may be solved explicitly in $\Gamma$-functions. 

Also note that the formula
$$
\Gamma(1+t) \, \Gamma(1-t) = \frac{2\pi i t}{e^{\pi i t} - e^{-\pi i t}}\,,
$$
and its relatives that were used above, show that the $\Gamma$-factors 
play the role of the Mukai's vector $\sqrt{\textup{Td} X}$ for the 
anti-linear inner product used in this paper. 
Independently, parallel $\Gamma$-factors arose in the work of Iritani \cite{Ir}.

\section{Connection Problem}

\subsection{} 

In this section, we solve the connection problem for the QDE between special points $q=0$ and $q=-1$, 
which may be called the DT and GW points, respectively. In plain English, the connection problem 
is to find the value of the matrix $\bY(q)$ at $q=-1$. We will see that the answer is 
given in terms of Macdonald polynomials. 

We will use $\bP^\lambda\in\cF\otimes\Q(q,t)$ to denote the monic Macdonald polynomial as 
defined, for example, in the book \cite{Mac}. Note that parameters $q$ and $t$, that are traditionally used 
for $\bP^\lambda$ are not our parameters $q$ and $t$. The matching of parameters will be 
discussed below. 

First, we explain the transformation that relates $\bP^\lambda$ to 
the polynomials $\bH^\lambda$ used by Haiman in his work on the $K$-theory of the Hilbert 
scheme, see \cite{Hai}. It is given by the formula 
\begin{equation}
  \label{Hform}
  \bH^\lambda(q,t) = t^{\bn(\lambda)}\, \prod_{\square\in\lambda} \left(1-q^{a(\square)} \, t^{-l(\square)-1}\right) \, \Upsilon \,\bP^\lambda(q,t^{-1})\,,
\end{equation}
where
$$
\Upsilon \lv \mu \rang  = \prod_\mu (1-t^{-\mu_i})^{-1} \, \lv \mu \rang \,. 
$$
and 
$$
\bn(\lambda) = \sum (i-1) \, \lambda_i \,.
$$
We relate Haiman's parameters $q$ and $t$ to our parameters by the identification 
$$
(q,t) = (T_1,T_2) \,.
$$
By their representation-theoretic meaning, they both define an element in the maximal 
torus of $GL(2)$ acting on $\C^2$ and its Hilbert schemes, so this identification is natural. 
Perhaps it would be even more natural to identify $q$ and $t$ with $T_i^{-1}$, but this 
amounts to a minor automorphism of symmetric functions.

\subsection{}

Let $\bH$ be the matrix with columns $\bH^\lambda$. We have the following 

\begin{Theorem}\label{Tconnect}
  \begin{equation}
    \label{connect}
     \bGm^{-1} \, \bY \, \bGd \Big|_{q=-1} = \frac{1}{(2\pi i)^{|\,\cdot\,|}} \, \bH 
  \end{equation}
\end{Theorem}

\noindent 
The proof of this theorem will take several steps.

\subsection{}

Let $H$ be an $\End(\cF)$-valued meromorphic function of $t_1$ and $t_2$ 
that satisfies 
\eqref{connect} in place of $\bH$. With this notation, $H=\bH$ is what we need to show. 

Formula \eqref{Sform} shows that $H$ is periodic in $t_i$ with period one, hence 
a meromorphic function of $T_1$ and $T_2$. Further, we claim that, in fact, $H$ is 
a rational function of the $T_i$'s. This follows from two 
following observations:
\begin{enumerate}
\item $H$ has finitely many poles;  
\item $H$ grows at most polynomially as $T_i\to0,\infty$\,.
\end{enumerate}
Indeed, the poles of $\bY^\lambda$ may occur only when 
$$
c(\lambda;t_1,t_2) = c(\mu;t_1,t_2)+n\,, \quad \mu\ne\lambda, 
\quad n=1,2,\dots\,.
$$
Since $c(\lambda;t_1,t_2)$ is a linear form with integer 
coefficients, these correspond to finitely many poles along divisors of the 
form 
$$
T_1^i T_2^j=1 \,.  
$$
In fact, these poles will be compensated by the zeros of $\bGd$, making $H$ 
a polynomial in the $T_i$'s, but we won't need this stronger result here. 

The growth of $H$ as $T_i\to0,\infty$ is estimated as in the proof of Theorem \ref{T_polynom}. 

\subsection{Asymptotics in the connection problem}

\subsubsection{}

Since $H$ is periodic, it may determined from studying its own asymptotics. Specifically, 
we will let $t_1\to+\infty$ keeping the level
$$
\kappa = t_1 + t_2
$$
fixed. In this limit, we have 
$$
 \bGm^{-1} \sim t_1^{\kappa|\,\cdot\,|-\ell(\,\cdot\,)} \, T_2^{-|\,\cdot\,|/2} \, \Upsilon \,. 
$$
Similarly, 
\begin{multline*}
\prod_{\textup{tangent weights $\bw$}} \frac{1}{\Gamma(\bw+1)} \sim (-2\pi i)^{|\lambda|}\, 
h_\lambda^{-1-\kappa}\, 
t_1^{-(1+\kappa)|\lambda|}\, \times \\  
T_1^{-\bn(\lambda')/2}\, 
T_2^{(\bn(\lambda)+|\lambda|)/2}\, 
\prod_{\square\in\lambda} \left(1-T_1^{a(\square)} T_2^{-l(\square)-1}\right)\,,
\end{multline*}
where $h_\lambda$ is the product of all hooklengths in the diagram of $\lambda$. 
Note that by Serre duality 
the tangent weights come in pairs $\bw_1,\bw_2$ such that $\bw_1+\bw_2=\kappa$. 

\subsubsection{}
Note that \eqref{branch} implies 
$$
q^{-c(\lambda,t_1,t_2)}\Big|_{q=-1} = T_1^{\bn(\lambda')/2} \, T_2^{\bn(\lambda)/2} 
$$
Together with the above asymptotics, this shows we should prove 
$$
(-1)^{|\lambda|} \, t_1^{-|\lambda|-\ell(\lambda)} \, h_\lambda^{-1-\kappa}\, \bY^\lambda(-1) \sim \bP^\lambda \,.
$$
To adjust for the $t_1$ scaling above, we introduce 
\begin{align*}
  \yb_\lambda &= (-1)^{|\lambda|}\, t_1^{-|\lambda|-\ell(\,\cdot\,)} \,
  \bY^\lambda \,, \\
\jb_\lambda &=  (-1)^{|\lambda|}\, t_1^{-|\lambda|-\ell(\,\cdot\,)}\,
  \bJ^\lambda \,, \\
\mb &= t_1^{-\ell(\,\cdot\,)} \, \MM_D \, t_1^{\ell(\,\cdot\,)}\,. 
\end{align*}
Note that from \eqref{dbJ} we have 
\begin{equation}\label{to_schur}
 \jb_\lambda \to 
h_\lambda\, s_\lambda\,,
\quad t_1\to \infty \,.
\end{equation}
The hook-length product in \eqref{to_schur} appears from the 
difference between the monic and the integral form of the Jack 
polynomial.

\subsubsection{}
At this point, we transformed the problem into showing that
\begin{equation}
  \label{asyyb}
  h_\lambda^{-1-\kappa}\, \yb_\lambda(-1) \sim \bP^\lambda \,, \quad t_1\to\infty \,.
\end{equation}
By its definition, the matrix $\bP$ is triangular in the basis $\{\jb_\mu\}$ 
with respect to the dominance order of partitions. 
Our next step is to show that the asymptotics of $\yb(-1)$ is similarly triangular.  

This will be done through a reformulation of the QDE in terms of an integral equation, following 
the strategy behind Levinson's theorem \cite{Lev}, see for example Section 1.4 in 
\cite{Easth}.  This technique is standard and is presented here mostly 
to make the material more accessible to algebraic geometers.

\subsubsection{}

It will be convenient to make a change of the independent variable via 
$$
q=-e^{-x} \,. 
$$
The new variable ranges from $0$ to $\infty$ as $q$ goes from $-1$ to $0$. 

We will write the differential equation 
$$
\frac{d}{dx} \, \Psi(x) = (-\mb-c(\lambda;t_1,t_2)) \, \Psi(x)
$$
satisfied by $\yb_\lambda$ in the form 
\begin{equation}
\frac{d}{dx} \, \Psi(x) = \left( D(x) +  R(x)\right) \, \Psi \,,
\label{dE}
\end{equation}
where the matrix $D(x)$ is diagonal in the basis of Jack polynomials 
$$
D(x) \, \jb_\mu = d_\mu (x) \, \jb_\mu \,,
$$
while the off-diagonal 
remainder matrix is linear in $\kappa$ 
and decays exponentially as 
$x\to +\infty$.

\subsubsection{}

Note that as $t_1\to\infty$ we have  
$$
d_\lambda(x;t_1,t_2) = t_1 \, \bff_2(\lambda) + O(1)\,, 
$$
where $\bff_2(\lambda)$ is the sum of contents of $\lambda$. 
Introduce a partial ordering $\lc$ on partitions by 
$$
\lambda \lc \mu \Leftrightarrow  \bff_2(\lambda) < \bff_2(\mu) \,. 
$$
This is a refinement of the dominance 
order on partitions, that is, 
$$
\lambda <  \mu \Rightarrow  \lambda \lc \mu  \,.
$$
We denote by $\Pi_{>}$ the projection onto $\jb_\mu$ with 
$\mu\gc\lambda$, that is 
$$
\Pi_{>} \, \jb_\mu = 
\begin{cases}
\jb_\mu\,, & \mu \gc \lambda\,,\\
0\,, & \textup{otherwise}\,.
\end{cases}
$$
We define the operator $\Pi_{<}$ similarly and set 
$$
\Pi_{\le} = 1 - \Pi_{>} \,, \quad  \Pi_{=} = 1 - \Pi_{>}-\Pi_{<} \,. 
$$
Given an operator $A$, we set 
$$
A_{>} = \Pi_{>} \, A \, \Pi_{>} \,.
$$
In English, this operator cuts out a corner of $A$ corresponding 
to rows and columns indexed partitions $\mu$ such that $\mu\gc \lambda$.

\subsubsection{}

Let the operator $A(x)$ be defined by
$$
A(x) \, \jb_\mu = e^{\int_0^x d_\mu(t) \, dt} \, \jb_\mu \,.
$$
It solves the equation 
$$
\frac{d}{dx} \, A(x) = D(x) \, A(x) \,.
$$

One constructs a particular 
solution of the differential equation \eqref{dE} 
as a solution of the integral equation 
\begin{multline}\label{iE}
\Psi(x) = e^{-\int_x^\infty d_\lambda(t)\, dt} \, 
\jb_\lambda 
\\
+ A(x)_{<} \int_0^x A(t)^{-1}\, R(t) 
\, \Psi(t)
\, dt \\ -
A(x)_{\le} \int_x^\infty A(t)^{-1} R(t)\, \Psi(t) \, dt \,. 
\end{multline}
Note that 
the convergence of the integral in the first line 
is insured by the exponential decay of $d_\lambda(x)$ as $x\to\infty$. 
Since 
$$
\frac{d}{dx} \, A_{<} = D A_{<}\,,
$$
it is immediate that a solution of \eqref{iE} solves \eqref{dE}.

\subsubsection{}

The fundamental role in the analysis of \eqref{iE} is played 
by the estimates 
\begin{align}
\left\| A(x)_{<} A(t)^{-1} \right \| &\le  
K_1(\kappa) \, e^{- t_1 (x-t)}\,, \notag \\
\left\| A(x)_{=} A(t)^{-1} \right \| &\le  
K_1(\kappa) \,, \label{estim} \\
\left\| A(x)_{>} A(t)^{-1} \right \| &\le  
K_1(\kappa) \, e^{- t_1 (t-x)}\,, \notag 
\end{align}
satisfied for all sufficiently large $t_1$. Here $K_1(\kappa)$ 
is constant depending on $\kappa$. In 
particular, it follows that a bounded solution of \eqref{iE} satisfies
$$
\Psi(x) \to \jb_\lambda\,, \quad x\to \infty \,.
$$

\subsubsection{}

Another crucial feature of the equation \eqref{iE} is 
\begin{equation}\label{Rze}
  \Pi_{=} \, R(x) \, \Pi_{=} = 0 \,.
\end{equation}
This is a geometric property. The matrix $R(x)$ is 
proportional to $\kappa=t_1+t_2$ and contains purely quantum 
parts of the quantum multiplication operator. The 
coefficient of $t_1+t_2$ is nonvanishing only 
if a chain of unbroken curves exists between two 
torus fixed points and the won't be any such 
chains if $\bff_2(\mu)=\bff_2(\lambda)$, see Section 3.8.2 of \cite{QCHS}.

\subsubsection{}

A bounded solution of \eqref{iE} may be constructed by successive
approximations 
\begin{multline*}
\Psi_{n+1}(x) = e^{-\int_x^\infty d_\lambda(t)\, dt} \, 
\jb_\lambda 
\\
+ A(x)_{<} \int_0^x A(t)^{-1}\, R(t) 
\, \Psi_n(t)
\, dt \\ -
A(x)_{\le} \int_x^\infty A(t)^{-1} R(t)\, \Psi_n(t) \, dt \,, 
\end{multline*}
starting with 
$$
\Psi_0(x)=e^{- \int_x^\infty d_\lambda(t)\, dt} \, 
\jb_\lambda \,.
$$
It follows from \eqref{estim} and \eqref{Rze} that 
\begin{align}
\left\| \Pi_{=} (\Psi_{n+1}-\Psi_{n}) \right\| & \le 
  K_1 K_2 \left\| \Pi_{\ne} (\Psi_{n}-\Psi_{n-1}) \right\| \,, 
\notag\\
\left\| \Pi_{\ne} (\Psi_{n+1}-\Psi_{n}) \right\| & \le 
\frac{ K_1 K_3}{t_1} \left\| \Psi_{n}-\Psi_{n-1} \right\|
\label{contract}
\end{align}
in the norm of $C([0,\infty))$, where 
$$
K_2(\kappa) = \left\|R(x)\right\|_{L^1([0,\infty))}\,, 
\quad 
K_3(\kappa) = \left\|R(x)\right\|_{C([0,\infty))} \,.  
$$ 
For $t_1\gg 0$ the iterations of \eqref{contract} are  
contracting, whence the convergence to a bounded solution 
$\Psi(x)$ of 
\eqref{iE}. Further, we have 
\begin{equation}
\Psi(x)\to \Psi_0(x) \,, \quad t_1\to + \infty \,. \label{asyPsi}
\end{equation}
%


\subsubsection{}

The asymptotics \eqref{asyPsi} implies that the connection matrix is 
diagonal in the $t_1\to\infty$ limit modulo the terms that are 
exponentially small as $x\to\infty$. In other words, this proves
that $\yb(-1)$ is asymptotically triangular in the basis $\{s_\mu\}$
with respect to the partial 
order $\gc$. 

Further, since $d_\mu(x)$ depends linearly on $\kappa$, 
the leading coefficients of $h_\lambda^{-1-\kappa}\, \yb_\lambda(-1)$ are of the form $K_4 K_5^\kappa$, where
$K_4$ and $K_5$ are real and positive. At the same time, they
should be periodic in $\kappa$, which yields $K_5=1$. 

Since a periodic function of $t_1$ is uniquely determined
by its asymptotics as $t_1\to+\infty$, this implies the corresponding
triangularity of $\Upsilon^{-1} H$. 

\subsection{}

We have 
$$
\bP^\lambda = s_\lambda + \dots 
$$
where dots stand for a linear combination of $s_\mu$ with $\lambda\gc \mu$. 
Therefore, 
$$
 H = \bH \, U 
$$
where the matrix $U$ is triangular with respect to $\gc$. Moreover, the 
diagonal entries of $U$ are positive real numbers that are independent 
of $t_1,t_2$. The unitarity of the connection and the formulas for the 
norm squares of the Macdonald polynomial imply
$$
U^* \, D \, U = D 
$$
where $U^* = \overline{U}^T$ and $D$ is the Gram matrix of the 
inner product \eqref{Herm_prod} in the basis $\bH^\lambda$. Since 
$D$ is diagonal, it follows that $U$ must be the identity matrix. 
This concludes the proof of 
Theorem \ref{Tconnect}.

\subsection{Scattering}\label{s_scatt}

Relation \eqref{q->1/q} implies the matrix
\begin{equation}\label{36gg}
(-1)^{\ell(\,\cdot\,)} \bY(q^{-1}) q^{c(\,\cdot\,)}
\end{equation}
is another fundamental solution of the quantum 
differential equation. It is natural to ask 
for the relation between \eqref{36gg} and the 
solution $\bY(q) q^{-c(\,\cdot\,)}$. Since, as 
$q\to 0,\infty$, our ODE becomes the Calogero-Sutherland
system, the 
$q\to q^{-1}$ transformation of the fundamental 
solution can be naturally interpreted as scattering
by the nonstationary terms.  

Since $q=-1$ is a fixed point of the involution $q\to q^{-1}$, 
we may use evaluation at $q=-1$ to compare the two 
solutions. {}From Theorem \ref{Tconnect}, we see that 
the scattering transformation essentially amounts to 
the action of the operator $(-1)^{\ell(\,\cdot\,)}$, i.e.\ the action of 
the standard symmetric function involution $\omega$,  in 
the basis $\{\bH^\lambda\}$. 

In particular, for integer levels $t_1+t_2 \in \Z$, we get 
the action of $\omega$ on Schur functions and so, up to 
normalization, scattering simply transposes the diagram 
$\lambda$ in that case.


\vspace{+10 pt}
\noindent
Department of Mathematics \\
Princeton University \\
Princeton, NJ 08544, USA\\
okounkov@math.princeton.edu \\

\vspace{+10 pt}
\noindent
Department of Mathematics\\
Princeton University\\
Princeton, NJ 08544, USA\\
rahulp@math.princeton.edu


\begin{thebibliography}{99}






\bibitem{BO} R.~Bezrukavnikov and A.~Okounkov, \emph{Monodromy of the 
QDE for the Hilbert scheme}, in preparation. 



\bibitem{jbrp} J.~Bryan and R.~Pandharipande, 
\emph{The local Gromov-Witten theory of curves}, math.AG/0411037.

\bibitem{Easth}
M.~S.~P.~Eastham,
\emph{The Asymptotic Solution of Linear Differential Systems. 
Applications of the Levinson Theorem.}, Clarendon Press, 
Oxford, 1989. 


\bibitem{Hai}
M.~Haiman, \emph{Combinatorics, symmetric functions and Hilbert schemes},
Current Developments in Mathematics, no. 1 (2002), 39-111.  


\bibitem{Ir}
H.~Iritani, 
\emph{An integral structure in quantum cohomology and mirror symmetry for toric orbifolds},
arXiv:0903.1463. 

\bibitem{Lev}
N.~Levinson,
\emph{The asymptotic nature of solutions of linear 
systems of differential equations}, Duke Math.\ J., 
\textbf{15} (1948) 111-126. 


\bibitem{LQW}
W.-P.~Li, Z.~Qin, W.~Wang, 
\emph{The cohomology rings of Hilbert schemes via Jack polynomials},
CRM Proceedings and Lecture Notes, vol.\ 38 (2004), 249--258. 



\bibitem{Mac}
I.~Macdonald, 
\emph{Symmetric functions and Hall polynomials}, 
The Clarendon Press, Oxford University Press, New York, 1995.


\bibitem{mnop} D.~Maulik, N.~Nekrasov, A.~Okounkov, and R.~Pandharipande, {\em Gromov-Witten
theory and Donaldson-Thomas theory I and II}, math.AG/0312059, math.AG/0406092. 



\bibitem{moop} D.~Maulik, A.~Oblomkov, A.~Okounkov, and R.~Pandharipande, {\em Gromov-Witten/Donaldson-Thomas correspondence for toric 3-folds}, arXiv:0809.3976.

\bibitem{QCHS} 
A.~Okounkov and R.~Pandharipande, 
\emph{Quantum cohomology of the Hilbert scheme of points in the plane}, 
arXiv:math/0411210. 

\bibitem{GWDT} 
A.~Okounkov and R.~Pandharipande, 
\emph{The local Donaldson-Thomas theory of curves}, arXiv:math/0512573. 


\bibitem{St}
R.~Stanley, 
\emph{Some combinatorial properties of Jack symmetric functions}, 
Adv.\ Math.\ \textbf{77} (1989), no.~1, 76--115. 


\bibitem{vass}
E.~Vasserot, 
\emph{Sur l'anneau de cohomologie du sch?ma de Hilbert de $\mathbf{C}\sp 2$},
C.~R.\ Acad.\ Sci.\ Paris S\'er. I Math.\ \textbf{332} (2001), no.~1, 7--12.

\end{thebibliography}
\end{document}